\documentclass[12pt,twoside]{article}
\usepackage{amsmath,amsbsy,amsfonts}
\usepackage[notref,notcite]{showkeys}
\newtheorem{theorem}{Theorem}[section]
\newtheorem{lemma}[theorem]{Lemma}
\newtheorem{proposition}[theorem]{Proposition}
\newtheorem{corollary}[theorem]{Corollary}
\newtheorem{remark}[theorem]{Remark}

\newcommand{\qed}{\hfill\rule{2mm}{2mm}}

\title{Nonlinear perturbations of a $p(x)$-Laplacian equation with critical growth in $\mathbb{R}^N$.}

\author{Claudianor O. Alves\footnote{C.O. Alves was partially supported by INCT-MAT, PROCAD, CNPq/Brazil
620150/2008-4 and 303080/2009-4, e-mail:coalves@dme.ufcg.edu.br}
\,\,\, and \,\,\,  Marcelo C. Ferreira \footnote{ e-mail:marcelo@dme.ufcg.edu.br}\\
Universidade Federal de Campina Grande\\
Unidade Acad\^emica de Matem\'atica e Estat\'istica\\
CEP:58429-900, Campina Grande - PB, Brazil.}

\date{}

\begin{document}
\maketitle



\begin{abstract}

We prove the existence of solution for a class of $p(x)$-Laplacian equations where the nonlinearity has a critical growth. Here, we consider two cases: the first case involves the situation where the variable exponents are periodic functions. The second one involves the case where the variable exponents are nonperiodic perturbations. 
\end{abstract}

\section{Introduction}

In this paper, we consider the existence of solution for the following class of equations 
$$ 
\left\{
\begin{array}{l}
-\Delta_{p(x)+\sigma(x)}u + (V(x)- W(x))| u |^{ p(x)+\sigma(x)-2 } u =  f(x,u)  \,\,\, \mbox{in} \ \mathbb{R}^{N},\\
u \geq 0, \,\, u\not= 0 \,\,\, \mbox{in} \,\,\, \mathbb{R}^{N} , \\
u \in W^{ 1,p(x)+\sigma(x) } \big( \mathbb{R}^{N} \big), 
\end{array}
\right.
\eqno{(P)}
$$
where $\Delta_{p(x)+\sigma(x)}$ is the $(p(x)+\sigma(x))$-Laplacian operator given by 
$$
\Delta_{p(x)+\sigma(x)}u= \text{div}(|\nabla u|^{p(x)+\sigma(x)-2} \nabla u),
$$
$f:\mathbb{R}^{N} \times \mathbb{R} \to \mathbb{R}$ is the function given by 
$$
f(x,t)=\mu |t|^{ q(x) - \tau(x) -2 } t + | t |^{ p^{ \ast }(x)-2 } t,
$$
$ \mu > 0 $ is a positive parameter, $ p, \sigma :  \mathbb{R}^{N} \to [0,+\infty)  $ are Lipschitz continuous functions and $ V,W, q, \tau  : \mathbb{R}^{N} \to [0, +\infty)$ are continuous functions verifying the following conditions: \\

\noindent The functions $p,q$ and $V$ are $\mathbb{Z}^{N}$-periodic, that is 
$$
p(x+y)=p(x), \, q(x+y)=q(x), \, V(x+y)=V(x) \,\,\, \forall x \in \mathbb{R}^{N} \,\,\, \mbox{and} \,\,\, \forall y \in \mathbb{Z}^{N}.
\eqno{(H_0)}
$$
Moreover, we assume also that 
$$
1 < p_- \leq p(x) \leq p_+ < N \,\,\, \forall x \in \mathbb{R}^{N}. \eqno{(H_1)}
$$
$$
p_+ < q_{-} \leq q(x) \ll p^{ \ast }(x) \,\,\, \forall x \in \mathbb{R}^{N}.  \eqno{(H_2)}
$$
\noindent There are $R > 0$ and $z \in \mathbb{R}^{N}$ with $B_R(z) \subset (0,1)^{N}$  such that
$$
\sigma(x), \tau(x) =0  \,\,\, \forall x \in \overline{B}^{c}_R(z).  \eqno{(H_3)}
$$
There are $m \in (1,N)$ and $R_1 > R$  with $B_{R_1}(z) \subset (0,1)^{N}$ such that 
$$
p(x)=m \,\,\, \forall x \in B_{R_1}(z). \eqno{(H_4)}
$$
$$
1 < (p+\sigma)_{-} \leq (p+\sigma)(x) \leq (p+\sigma)_+ < N \,\,\, \forall x \in \mathbb{R}^{N}. \eqno{(H_5)}
$$
$$
(p+\sigma)_+ < (q-\tau)_{-} \leq (q-\tau)(x) \ll (p+\sigma)^{ \ast }(x) \,\,\, \forall x \in \mathbb{R}^{N}.  \eqno{(H_6)}
$$
$$
W(x) \to 0 \,\,\, \mbox{as} \,\,\, |x| \to +\infty. \eqno{(W_0)}
$$
$$
\inf_{ x \in \mathbb{R}^{N} } V(x) = V_0>0.  \eqno{(V_0)}
$$ 
and
$$
\inf_{ x \in \mathbb{R}^{N} }(V(x)-W(x)) = U_0>0.  \eqno{(WV_0)}
$$

Here, the notation $h  \ll  g$  means that  $\displaystyle \inf_{ x \in \mathbb{R}^{N}}(g(x)-h(x))>0$, \linebreak $h_- = \displaystyle \mbox{ess}\inf_{x \in \mathbb{R}^{N}}h(x)$, $h_+ = \displaystyle \mbox{ess}\sup_{x \in \mathbb{R}^{N}}h(x)$ and $h^{*}(x)=\frac{Nh(x)}{N-h(x)}$.

\vspace{0.5 cm}

The study of problems with variable exponents has received a special attention at the last years, because this class of problems appears in various mathematical models, such as, \\

 \noindent $\bullet$ {\it Electrorheological  fluids:} see Acerbi \& Mingione \cite{Acerbi1,Acerbi2}, Antontsev \& Rodrigues \cite{Antontsev1} and Ruzicka \cite{Ru}, \\

 \noindent $\bullet$ {\it Nonlinear Darcy law in porous medium:} see Antontsev \& Shmarev \cite{Antontsev2,Antontsev3} \\

 \noindent $\bullet$ {\it Image Processing:} see Chambolle \& Lions \cite{CLions} and Chen, Levine \& Rao \cite{Chen}.

\vspace{0.5 cm}

Motivated by the presence of variable exponents in the applications
above and after Kovacik \& R\'akosn\'ik have shown  some properties of the spaces $L^{p(x)}$
and $W^{1,p(x)}$ in \cite{KR}, a lot of research has been done concerning
these kinds of problems, see, for example, Alves \cite{Alves1,Alves2}, Alves \& Souto \cite{AlvesSouto}, Fan \cite{Fan}, Fu \& Zhang \cite{FuZa}, Fan \& Han \cite{FanHan},  Krist\'aly, Radulescu \& Varga \cite{KRV},  Mihailescu \& Radulescu \cite{MR} and references therein.

In \cite{ACM}, Alves, Carri\~ao \& Miyagaki have considered the existence of solution for problem $(P)$, for a case where the exponents $p$ and $q$ are constants and $\sigma=\tau=0$. More precisely, in \cite{ACM}, the following problem was studied 
$$ 
\left\{
\begin{array}{l}
-\Delta u +(V(x)-W(x))u=\mu u^{q-1} + u^{2^{*}-1} \,\,\, \mbox{in} \,\,\, \mathbb{R}^{N}, \\
u >0 \,\,\, \mbox{in} \,\,\, \mathbb{R}^{N}, \\
u \in H^{1}(\mathbb{R}^{N}), 
\end{array}
\right. \eqno{(P_1)}
$$
where $\mu>0$ is a positive parameter, $ q \in (2,2^{*})$ and $W \in L^{\frac{N}{2}}(\mathbb{R}^{N})$.
In that paper, the authors used variational methods combined with a well known result due to Lions \cite{Lions} and showed the existence of solutions for all $\mu > 0$. An important point in that work is the fact that the continuous embedding $W^{1,p}(\mathbb{R}^{N}) \hookrightarrow L^{p^*}(\mathbb{R}^{N})$ has a best Sobolev constant, denoted by S, which is assumed by a special class of functions. When the function $p$ is not constant, we do not have this information, and thus new arguments and estimates are necessary.

In \cite{Fan}, Fan has considered a class of nonperiodic perturbations like problem $(P)$, however in that paper the nonlinearity has a subcritical growth. More precisely, the problem studied was the following 
$$ 
\left\{
\begin{array}{l}
-div(|\nabla u|^{p(x) + \sigma(x)-2}\nabla u)+ a(x) |u|^{p(x) + \sigma(x)-2}u= |u|^{ -\tau(x) }f(x,u)\,\,\, \mbox{in} \,\,\, \mathbb{R}^{N}, \\
u \geq 0, \,\, u\not= 0 \,\,\, \mbox{in} \,\,\, \mathbb{R}^{N}, \\
u \in W^{1,p(x)}(\mathbb{R}^{N}), 
\end{array}
\right.  \eqno{(P_2)}
$$

\vspace{0.5 cm}
\noindent where $p,q,a,\sigma,\tau: \mathbb{R}^{N} \to \mathbb{R}$ are continuous functions with $p$ and $q$ being $\mathbb{Z}^{N}$-periodic functions and $f:\mathbb{R}^{N} \times \mathbb{R} \to \mathbb{R}$ being a continuous function with subcritical growth.  The main tool used was the variational method, more precisely, some characterizations of the mountain pass corresponding to the energy functional associated with problem $(P_2)$.

\vspace{0.5 cm}

Motivated by papers \cite{ACM} and  \cite{Fan}, we will show the existence of solution for problem $(P)$ by using the variational method. Here, we look for critical points of the energy functional associated with $(P)$ given by 
$$
I(u)=\int_{\mathbb{R}^{N}} \frac{1}{(p(x)+\sigma(x))} \left( | \nabla u |^{ p(x)+\sigma(x) } + (V(x)-W(x)) | u |^{ p(x)+\sigma(x) } \right) - \Psi(u) - J(u) 
$$
where
$$
\Psi(u)= \mu \int_{\mathbb{R}^{N}} \frac{1}{(q(x)-\tau(x))} | u |^{ q(x)-\tau(x) } \,\,\, \mbox{and} \,\,\, J(u)=\int_{\mathbb{R}^{N}} \frac{1}{p^{ \ast }(x)} |u|^{ p^{ \ast}(x) },  
$$
for all  $u \in W^{ 1,p(x)+\sigma(x) }( \mathbb{R}^{N} )$.  

Hereafter, we will consider the following norm in $ W^{ 1,p(x)+\sigma(x) } \big( \mathbb{R}^{N} \big) $ 
$$ 
      \| u \| = \inf \left\{ \alpha > 0; \, \rho \big( \alpha^{ -1 }u \big) \le 1 \right\},
$$
where 
$$
      \rho( u ) = \int_{\mathbb{R}^{N}} (| \nabla u |^{ p(x)+\sigma(x) }+ (V(x)-W(x)) | u |^{ p(x)+\sigma(x) }).
$$

Using well known arguments, we have that $ I \in C^1 \Big( W^{ 1,p(x)+\sigma(x) } \big( \mathbb{R}^{N} \big), \mathbb{R} \Big) $ with 
$$ 
I'( u ) v  = \int_{\mathbb{R}^{N}}( | \nabla u |^{ p(x)+\sigma(x)-2 } \nabla u \nabla v + (V(x)-W(x)) | u |^{ p(x)+\sigma(x)-2 } u v ) - \Psi'(u)v- J'(u)v,  
$$
where  
$$
\Psi'(u)v= \mu \int_{\mathbb{R}^{N}} | u |^{ q(x)-\tau(x)-2 } u v \,\,\, \mbox{and} \,\,\, 
J'(u)v= \int_{\mathbb{R}^{N}} | u |^{ p^{ \ast }(x)-2 } u v,
$$
for all $ u,v \in W^{ 1,p(x)+\sigma(x)}( \mathbb{R}^{N} ) $.

\vspace{0.5 cm}

Our main result is the following:

\begin{theorem} \label{T1} Assume $(H_0)-(H_6),(V_0),(W_0)$ and $ (WV_0) $. Then, there is $\mu^*>0$, such that problem $(P)$ has a nonnegative ground-state solution for all $\mu \geq \mu^{*}$. 
\end{theorem}

In what follows, we mean that a solution $u$ of $(P)$ is a {\it  ground-state solution}, if it is a least energy solution, that is, if for any nontrivial solution $v$ of  $(P)$, we have that $I(u) \leq I(v)$. 

\vspace{0.5 cm}

\noindent {\bf Notation:} The following notations will be used in the present work: \\

\noindent $\bullet$ \,\, $C$ and $C_i$ will denote generic positive constant, which may vary from line to line.  \\

\noindent $\bullet$ \,\, In all the integrals we omit the symbol $dx$.

\section{The periodic problem}

In this section, we study the existence of a ground-state solution for the periodic problem related to $(P)$ given by
$$ 
\left\{
\begin{array}{l}
-\Delta_{p(x)}u + V(x)| u |^{ p(x)-2 } u = \mu |u|^{ q(x)-2 } u + | u |^{ p^{ \ast }(x)-2 } u  \,\,\, \mbox{in} \ \mathbb{R}^{N},\\
\mbox{}\\
u \in W^{ 1,p(x)} \big( \mathbb{R}^{N} \big), u \neq 0 .
\end{array}
\right.
\eqno{(P_\infty)}
$$

Our goal in this section is to prove the following result

\begin{theorem} \label{T2} Assume $(H_0)-(H_2)$ and $(V_0)$. Then, there is $\mu_{\infty}>0$, such that problem $(P_\infty)$ has a nonnegative ground-state solution for all $\mu \geq \mu_{\infty}$.
\end{theorem}

The energy functional $I_{\infty}: W^{ 1,p(x)} \big( \mathbb{R}^{N} \big) \to \mathbb{R}$ associated with $(P_\infty)$ is given by 
$$
I_\infty(u)=\int_{\mathbb{R}^{N}} \frac{1}{p(x)} \left( | \nabla u |^{ p(x) } + V(x) | u |^{ p(x) } \right) - \mu \int_{\mathbb{R}^{N}} \frac{1}{q(x)} | u |^{ q(x)} - \int_{\mathbb{R}^{N}} \frac{1}{p^{ \ast }(x)} |u|^{ p^{ \ast}(x) }.
$$
A direct computation shows that $I_\infty \in C^1(W^{1,p(x)}(\mathbb R^N),\mathbb{R})$ with 
$$ 
\begin{array}{l}
I_\infty'( u ) v  = \displaystyle \int_{\mathbb{R}^{N}}( | \nabla u |^{ p(x)-2 } \nabla u \nabla v + V(x) | u |^{ p(x)-2 } u v ) - \mu \int_{\mathbb{R}^{N}} | u |^{ q(x)-2 } u v \\ 
\mbox{}\\
\hspace{2 cm} - \displaystyle \int_{\mathbb{R}^{N}} | u |^{ p^{ \ast }(x)-2 } u v,  \\
\end{array}
$$
for all $ u,v \in W^{ 1,p(x)}( \mathbb{R}^{N} ) $.

Using standard arguments, it is easy to prove that $ I_{ \infty } $ satisfies the mountain pass geometry, this way, there exists $ (u_n) \subset W^{ 1,p(x) } \big( \mathbb{R}^{N} \big)  $ verifying 
\begin{equation} \label{NIVELZ}
      I_{ \infty }( u_n ) \to c_{ \infty } \qquad \text{and} \qquad I'_{ \infty }( u_n ) \to 0,
\end{equation}
where 
$$ 
      c_{ \infty} = \inf_{ \gamma \in \Gamma } \max_{ t \in [0,1] } I_{ \infty }( \gamma(t) ),
$$
and
$$ 
\Gamma = \left\{ \gamma \in C \left( [0,1], W^{ 1,p(x) } \big( \mathbb{R}^{N} \big) \right); \, \gamma(0) = 0 \ \text{and} \ I_{ \infty }( \gamma(1) ) \le 0 \right\}. 
$$
The level $c_\infty$ is called the mountain pass level of the functional $I_{\infty}$. An important point that we would like to mention is the fact that 
\begin{equation} \label{LIMITEDONIVEL}
     c_{ \infty } \to 0 \,\,\, \mbox{as} \,\,\, \mu \to +\infty  \,\,\, ( \mbox{see} \,\,\, \cite{Alves1} ).
\end{equation} 
Using the above information, in the present paper, we fix $ \mu_{ \infty } > 0 $ such that  
\begin{equation} \label{ESTIMATIVAPASSO}
      c_{ \infty } < \min \left \{ \theta \left( \frac{1}{K} \right)^{ \frac{1}{\theta} }, \frac{1}{2K^{ p_+ }} \nu \right \} \,\,\, \forall \mu \ge \mu_{ \infty },
\end{equation} 
where 
\begin{equation} \label{theta}
\theta = 1/p_+ - 1/p^{ \ast}_-, \ \nu = 1/p_+ - 1/q_-,  
\end{equation}
and $K \ge 1$ is fixed satisfying  
$$  
| u |_{ p^{ \ast }(x) } \le K \| u \|, \,\,\,\,  \forall u \in W^{ 1,p(x) } \big( \mathbb{R}^{N} \big).
$$

\vspace{0.5 cm}

Next, we will make a brief review about the spaces $L^{p(x)}(\mathbb{R}^{N})$ and $W^{1,p(x)}(\mathbb{R}^{N})$.

\subsection{Variable exponent Lebesgue and Sobolev spaces}

In this subsection, we recall some results on variable exponent Lebesgue and Sobolev spaces found in \cite{FZ,FSZ} and their references. 

Let $h\in L^{\infty}(\mathbb{R}^{N})$ with $ h_-> 1$. The variable exponent Lebesgue space
$L^{h(x)}(\mathbb{R}^{N})$ is defined by
\[
L^{h(x)}(\mathbb{R}^{N})=\left\{
u:\mathbb{R}^{N} \to  \mathbb{R} \left\vert \,u\text{ is
measurable and }\int_{\mathbb{R}^{N}}\left\vert u\right\vert
^{h(x)}<\infty\right.  \right\}
\]
endowed with the norm
\[
\left\vert u\right\vert _{h(x)}=\inf\left\{  \lambda>0\left\vert
\,\int_{\mathbb{R}^{N}}\left\vert \frac{u}{\lambda}\right\vert
^{h(x)} \leq1\right.  \right\}  \text{.}
\]
The variable exponent Sobolev space is defined by 
\[
W^{1,h(x)}(\mathbb{R}^{N})=\left\{  u\in
L^{h(x)}(\mathbb{R}^{N})\left\vert \,\left\vert \nabla u\right\vert
\in L^{h(x)}(\mathbb{R}^{N})\right. \right\}
\]
with the norm
\[
\left\Vert u\right\Vert _{1,h(x)}=\left\vert u\right\vert
_{h(x)}+\left\vert \nabla u\right\vert _{h(x)}\text{.}
\]
When $M \in L^{\infty}(\mathbb{R}^{N})$ and $ M_- > 0 $, the norm
\begin{equation}
\left\Vert u\right\Vert =\inf\left\{  \lambda>0\left\vert \,\int
_{\mathbb{R}^{N}}\left(  \left\vert \frac{\nabla
u}{\lambda}\right\vert ^{h(x)}+M(x)\left\vert
\frac{u}{\lambda}\right\vert ^{h(x)}\right) \leq1\right.
\right\}  \label{n}
\end{equation}
is equivalent to norm $\left\Vert \, \cdot \,\right\Vert _{1,h(x)}$.
With these norms, the spaces $L^{h(x)}(\mathbb{R}^{N})$ and
$W^{1,h(x)}(\mathbb{R}^{N})$ are reflexive and separable Banach spaces.

\begin{proposition}
\label{p1}The functional $\zeta:W^{1,h(x)}(\mathbb{R}^{N}) \to \mathbb{R}$ defined by
\begin{equation}
\zeta(u)=\int_{\mathbb{R}^{N}}\left(  \left\vert \nabla u\right\vert
^{h(x)}+M(x)\left\vert u\right\vert ^{h(x)}\right)
\text{,} \label{ps}
\end{equation}
has the following properties:

\begin{enumerate}
\item[\emph{(i)}] If $\left\Vert u\right\Vert \geq1$, then $\left\Vert u\right\Vert
^{h_{-}}\leq\zeta(u)\leq\left\Vert u\right\Vert ^{h_+}$.

\item[\emph{(ii)}] If $\left\Vert u\right\Vert \leq1$, then $\left\Vert u\right\Vert
^{h_+}\leq\zeta(u)\leq\left\Vert u\right\Vert ^{h_{-}}$.
\end{enumerate}
In particular, $\zeta(u)=1$ if, and only if, $\left\Vert u \right\Vert =1$ and, for $ (u_n) \subset W^{ 1,h(x) }( \mathbb R^N ) $, $\left\Vert u_{n}\right\Vert
\rightarrow0$ if, and only if, $\zeta( u_{n}) \rightarrow0$.
\end{proposition}

\begin{remark}
\label{r1}For the functional
$\xi:L^{h(x)}(\mathbb{R}^{N})\rightarrow \mathbb{R}$ given by 
\[
\xi(u)=\int_{\mathbb{R}^{N}}\left\vert u\right\vert
^{h(x)}\text{,}
\]
the conclusion of Proposition \ref{p1} also holds, for example, if $ (u_n) \subset L^{h(x) }( \mathbb R^N ) $,  $\left\vert u_{n}\right\vert
_{h(x)}\rightarrow0$ if, and only if $\xi(u_{n})\rightarrow0$. Moreover, from $(i)$ and $(ii)$, 
\begin{equation}
\left\vert u\right\vert _{h(x)}\leq \max\left\{
\left(\int_{\mathbb{R}^{N}}\left\vert u\right\vert
^{h(x)}\right)  ^{1/h_{-}},\left(  \int_{\mathbb{R}^{N}
}\left\vert u\right\vert ^{h(x)}\right)
^{1/h_+} \right\} \text{.} \label{in}
\end{equation}

\end{remark}

Related to the Lebesgue space $L^{h(x)}(\mathbb{R}^{N})$, we have the following generalized H\"{o}lder's inequality.

\begin{proposition}
[{\cite[p.9]{Mu}}]\label{h} For $h\in L^{\infty}( \mathbb{R}^{N}) $ with $ h_- >1$, let \linebreak $h^{\prime }:\mathbb{R}^{N} \to \mathbb{R}$ be such that
\[
\frac{1}{h( x) }+\frac{1}{h^{\prime}( x) }=1\text{,\qquad a.e. }x\in
\mathbb{R}^{N}\text{.}
\]
Then, for any $u\in L^{h( x) }( \mathbb{R}^{N})
$ and  $v\in L^{h^{\prime}( x) }( \mathbb{R}^{N}) $,
\begin{equation}
\left\vert \int_{\mathbb{R}^{N}}uv\,\right\vert
\leq \left( \frac{1}{h_-} + \frac{1}{h'_-} \right)\left\vert u\right\vert _{h( x) }\left\vert v\right\vert
_{h^{\prime}( x) }\text{.} \label{hi}
\end{equation}
\end{proposition}

\begin{proposition}
[{\cite[Theorems 1.1, 1.3]{FSZ}}]\label{pp1}Let $h:\mathbb{R}^{N}
\to \mathbb{R}$ be a Lipschitz continuous satisfying $ 1 < h_- \leq h_+ < N $ and $t:\mathbb{R}^{N} \rightarrow \mathbb{R}$ be a measurable function.

\begin{enumerate}
\item[\emph{(i)}] If $h\leq t\leq h^{\ast}$, the embedding
$W^{1,h(x)}(\mathbb{R}^{N})\hookrightarrow
L^{t(x)}(\mathbb{R}^{N})$ is continuous.

\item[\emph{(ii)}] If $h\leq t\ll h^{\ast}$, the embedding
$W^{1,h(x)}(\mathbb{R}^{N})\hookrightarrow L_{\mathrm{loc}}^{t(x)}
(\mathbb{R}^{N})$  is compact.
\end{enumerate}
\end{proposition}

\subsection{Preliminary results}

\begin{lemma} \label{aew convergence}
   Let $ ( v_n ) $ be a $ (PS)_d $ sequence for $ I_{ \infty } $ and $ v \in W^{ 1,p(x) } \big( \mathbb{R}^{N} \big) $ such that 
   $$
      v_n \rightharpoonup v.
   $$
   Then,
   $$ 
      I'_{ \infty }(v) = 0.
   $$
Thus, if $ v \ne 0 $, $ v $ is a nontrivial solution of $ ( P_{ \infty } ) $. 
\end{lemma}

\noindent {\bf Proof.} 
Following a standard reasoning, it is sufficient to show that, up to a subsequence ,
$$
\nabla v_n( x ) \to \nabla v( x ) \,\,\, \mbox{a.e in} \,\,\,  \ \mathbb{R}^{N}.
$$
We begin observing that, up to a subsequence, there exist nonnegative measures $ \mathfrak{m} $ and $ \mathfrak{n} $ in $ \mathcal{ M }\big( \mathbb{R}^{N} \big) $ such that
\begin{equation} \label{MED1}
| \nabla v_n |^{ p(x) }  \rightharpoonup \mathfrak{m} \ \text{in} \ \mathcal{M} \big( \mathbb{R}^{N} \big) 
\end{equation}
and
\begin{equation} \label{MED2}
| v_n |^{ p^{ \ast }(x) }  \rightharpoonup \mathfrak{n} \ \text{in} \ \mathcal{M} \big( \mathbb{R}^{N} \big).
\end{equation}
By using a concentration compactness principle found in \cite{FuZa}, there exists a countable index set $ \mathfrak{I} $ such that  
   \begin{gather*}
      \mathfrak{n} = | v |^{ p^{ \ast }(x) } \, dx + \sum_{ i \in \mathfrak{I} } \mathfrak{n}_i \delta_{ x_i }, \\
      \mathfrak{m} \ge | \nabla v |^{ p(x) } \, dx + \sum_{ i \in \mathfrak{I} } \mathfrak{m}_i \delta_{ x_i }, 
 \end{gather*}      
and
$$
     \mathfrak{n}_i \le S \max \left\{ \mathfrak{m}_i^{ \frac{p^{ \ast}_+}{p_-} }, \mathfrak{m}_i^{ \frac{p^{ \ast}_-}{p_+} } \right\}
$$  
where $ ( \mathfrak{n}_i )_{ i \in \mathfrak{I}}, ( \mathfrak{m}_i )_{ i \in \mathfrak{I} } \subset [0, \infty) $ and $ ( x_i )_{ i \in \mathfrak{I} } \subset \mathbb{R}^{N}. $ The constant $ S $ is given by
$$
      S = \sup_{ u \in W^{ 1,p(x) } ( \mathbb{R}^{N} ) \atop \| u \| \le 1 } \int_{\mathbb{R}^{N}} |u|^{ p^{ \ast }(x) }.
$$
Our first task is to prove that 
$$ 
\mathfrak{m_i} = \mathfrak{n_i}, \, \forall i \in \mathfrak{I}.
$$
For this, let $ \varphi \in C^{ \infty }_0 \big( \mathbb{R}^{N} \big) $ such that 
$$
\varphi(x) = 1 \,\,\,  \text{in} \ B_1( 0 ), \ \varphi(x) = 0  \ \text{in} \  B_2^c( 0 ) \ \text{and} \ 0 \le \varphi(x) \le 1 \, \forall \, x \in \mathbb{R}^{N}.
$$  
Fixed $ i \in \mathfrak{I} $, we consider for each $ \epsilon > 0 $  
$$
\varphi_{ \epsilon } (x) = \varphi \left( \frac{x - x_i}{\epsilon} \right) \,\,\, \forall x \in \mathbb{R}^{N}.
$$
Since $ ( v_n ) $ is bounded in $ W^{ 1,p(x) } \big( \mathbb{R}^{N} \big) $, the sequence $ ( \varphi_{ \epsilon } v_n ) $ is also bounded in $ W^{ 1,p(x) } \big( \mathbb{R}^{N} \big) $. Thus,
$$
      I'_{ \infty }(v_n)( \varphi_{ \epsilon } v_n ) = o_n(1),
$$
that is, 
   \begin{multline*}
      \int_{\mathbb{R}^{N}} (\varphi_{ \epsilon } | \nabla v_n |^{ p(x) } + v_n | \nabla v_n |^{ p(x)-2 } \nabla v_n  \nabla \varphi_{ \epsilon }) + \int_{\mathbb{R}^{N}} V(x) | v_n |^{ p(x) } \varphi_{ \epsilon }  \\
      = \mu \int_{\mathbb{R}^{N}} | v_n |^{ q(x) } \varphi_{ \epsilon }  + \int_{\mathbb{R}^{N}} | v_n |^{ p^{ \ast }(x) } \varphi_{ \epsilon } + o_n(1).
   \end{multline*}
   Taking the limit as $ n \to \infty $, the weak convergence  of $ (| \nabla v_n |^{ p(x) } ) $ and $ (| v_n |^{ p^{ \ast } (x) } )$ in $\mathcal {M}(\mathbb{R}^N)$ combined with the Lebesgue dominated convergence theorem leads to 
   \begin{multline} \label{epsilon limit} 
      \int_{\mathbb{R}^{N}} \varphi_{ \epsilon } \, d \mathfrak{m} + \limsup_n \int_{\mathbb{R}^{N}} v_n | \nabla v_n |^{ p(x)-2 } \nabla v_n \nabla \varphi_{ \epsilon }  + \int_{\mathbb{R}^{N}} V(x) | v |^{ p(x) } \varphi_{ \epsilon }  \\
      = \mu \int_{\mathbb{R}^{N}} | v |^{ q(x) } \varphi_{ \epsilon }+ \int_{\mathbb{R}^{N}} \varphi_{ \epsilon } \, d \mathfrak{n}.
   \end{multline}
Using H\"older's inequality and the boundedness of $(v_n)$ in $W^{1,p(x)}(\mathbb{R}^{N})$, we get
   \begin{align*}
      & \left| \int_{\mathbb{R}^{N}} v_n | \nabla v_n |^{ p(x)-2 } \nabla v_n \cdot \nabla \varphi_{ \epsilon }  \right| \\
      & \phantom{ \left| \right. } \le \int_{\mathbb{R}^{N}} \left| \nabla v_n \right|^{ p(x)-1 } \left| v_n \nabla \varphi_{ \epsilon } \right|  
        \le C \left| \left| \nabla v_n \right|^{ p(x)-1 } \right|_{ p'(x) } \big| v_n \left| \nabla \varphi_{ \epsilon } \right| \big|_{ p(x) } \\
      & \phantom{ \left| \right. } \le C \max \left\{ \left( \int_{\mathbb{R}^{N}} | v_n |^{ p(x) } \left| \nabla \varphi_{ \epsilon } \right|^{ p(x) }  \right)^{ \frac{1}{p_-} }, \left( \int_{\mathbb{R}^{N}} | v_n |^{ p(x) } \left| \nabla \varphi_{ \epsilon } \right|^{ p(x) } \right)^{ \frac{1}{p_+} } \right\}, 
   \end{align*}
where $ p'(x) = \frac{ p(x) }{ p(x)-1 }  \,\,\, \forall x \in \mathbb{R}^{N} $. Therefore, by Sobolev embedding 
   \begin{multline*}
      \limsup_n \left|\int_{\mathbb{R}^{N}} v_n | \nabla v_n |^{ p(x)-2 } \nabla v_n \cdot \nabla \varphi_{ \epsilon }  \right| \\
        \le C \max \left\{ \left( \int_{\mathbb{R}^{N}} | v |^{ p(x) } \left| \nabla \varphi_{ \epsilon } \right|^{ p(x) } \right)^{ \frac{1}{p_-} }, \left( \int_{\mathbb{R}^{N}} | v |^{ p(x) } \left| \nabla \varphi_{ \epsilon } \right|^{ p(x) }  \right)^{ \frac{1}{p_+} } \right\}. 
   \end{multline*}
Furthermore, by H\"older's inequality
   $$
      \int_{\mathbb{R}^{N}} | v |^{ p(x) } \left| \nabla \varphi_{ \epsilon } \right|^{ p(x) }
        \le C \left| | v |^{ p(x) } \right|_{ L^{ \frac{N}{N-p(x)} } \big( B_{ 2 \epsilon }(x_i) \big) } \big| \left| \nabla \varphi_{ \epsilon } \right|^{ p(x) } \big|_{ L^{ \frac{N}{p(x)} } \big( B_{ 2 \epsilon }(x_i) \big) }. \\
   $$
   Once that
   $$
      \int_{ B_{ 2 \epsilon }(x_i) } \left| \nabla \varphi_{ \epsilon } \right|^N  = \int_{ B_2(0) } \left| \nabla \varphi \right|^N ,
   $$
   we derive 
   \begin{multline*}
      \big| \left| \nabla \varphi_{ \epsilon } \right|^{ p(x) } \big|_{ L^{ \frac{N}{p(x)} } \big( B_{ 2 \epsilon }(x_i) \big) } \\
      \phantom{ \big| } \le \max \left\{ \left( \int_{ B_{ 2 \epsilon }(x_i) }  \left| \nabla \varphi_{ \epsilon } \right|^N  \right)^{ \frac{1}{\left( \frac{N}{p} \right)_-}  }, \left( \int_{ B_{ 2 \epsilon }(x_i) }  \left| \nabla \varphi_{ \epsilon } \right|^N  \right)^{ \frac{1}{\left( \frac{N}{p} \right)_+}  } \right\} \leq C
   \end{multline*}
 for some positive constant $C$, which is independent of $\epsilon$.  Thereby, 
   $$
      \int_{\mathbb{R}^{N}} | v |^{ p(x) } \left| \nabla \varphi_{ \epsilon } \right|^{ p(x) } \le C \left| | v |^{ p(x) } \right|_{ L^{ \frac{N}{N-p(x)} } \big( B_{ 2 \epsilon }(x_i) \big) },
   $$
   and so
   \begin{multline*}
     \limsup_n \left|\int_{\mathbb{R}^{N}} v_n | \nabla v_n |^{ p(x)-2 } \nabla v_n \cdot \nabla \varphi_{ \epsilon } \right| \\
      \le C \max \left\{ \left| | v |^{ p(x) } \right|_{ L^{ \frac{N}{N-p(x)} } \big( B_{ 2 \epsilon }(x_i) \big) }^{ \frac{1}{p_-} }, \left| | v |^{ p(x) } \right|_{ L^{ \frac{N}{N-p(x)} } \big( B_{ 2 \epsilon }(x_i) \big) }^{ \frac{1}{p_+} } \right\}. \\
   \end{multline*}
   But,
   \begin{multline*}
      \left| | v |^{ p(x) } \right|_{ L^{ \frac{N}{N-p(x)} } \big( B_{ 2 \epsilon }(x_i) \big) } \\
      \le \max \left\{ \left( \int_{ B_{ 2 \epsilon }(x_i) }  | v |^{ p^{ \ast }(x) } \right)^{ \frac{1}{\left( \frac{N}{N-p} \right)_- } }, \left( \int_{ B_{ 2 \epsilon }(x_i) }  | v |^{ p^{ \ast }(x) }  \right)^{ \frac{1}{\left( \frac{N}{N-p} \right)_+ } } \right\}
   \end{multline*} 
from where it follows that  
   $$
      \lim_{\epsilon \to 0}\limsup_n \left| \int_{\mathbb{R}^{N}} v_n | \nabla v_n |^{ p(x)-2 } \nabla v_n \nabla \varphi_{ \epsilon } \right| = 0
   $$
implying that
$$
      \lim_{\epsilon \to 0}\limsup_n  \int_{\mathbb{R}^{N}} v_n | \nabla v_n |^{ p(x)-2 } \nabla v_n \nabla \varphi_{ \epsilon } = 0.
 $$
Now, taking the limit as $ \epsilon \to 0 $ in  \eqref{epsilon limit},  we get
   \begin{equation} \label{m_i = n_i}
      \mathfrak{m}_i = \mathfrak{m}( x_i ) = \mathfrak{n}( x_i ) = \mathfrak{n}_i.
   \end{equation}
Once that 
$$   
\frac{ p^{\ast}_- }{ p_+ } \le \frac{ p^{ \ast }_+ }{ p_- }, 
$$
we have that 
\begin{equation} \label{eta1}
      \mathfrak{n}_i^{ \frac{p_+}{p^{ \ast }_-} } \le \left( S^{ \frac{p_+}{p^{ \ast }_-} } + S^{ \frac{p_-}{{p^{ \ast }_+}} } \right) \mathfrak{m}_i, \ \text{if} \,\,\, \mathfrak{m}_i <1
 \end{equation}
and
\begin{equation} \label{eta2}
\mathfrak{n}_i^{ \frac{p_-}{p^{ \ast }_+} } \le \left( S^{ \frac{p_+}{p^{ \ast }_-}} + S^{ \frac{p_-}{{p^{ \ast }_+}} } \right) \mathfrak{m}_i \,\,\, \mbox{if} \,\,\, \mathfrak{m}_i \geq 1.
\end{equation}
Thus, from \eqref{m_i = n_i} - \eqref{eta2}, if $ \mathfrak{n}_i > 0 $ for some $ i \in \mathfrak{I} $, there exists $\alpha>0$, which is independent of $ i $, such that
\begin{equation} \label{eta3}
\mathfrak{n}_i \geq \alpha.
\end{equation}
Recalling that
  \begin{equation} \label{n_i sum}
     \sum_{ i \in \mathfrak{I} \atop \mathfrak{m}_i < 1 } \mathfrak{n}_i^{ \frac{ p_+ }{ p^{\ast}_- } } + 
       \sum_{ i \in \mathfrak{I} \atop \mathfrak{m}_i \ge 1 } \mathfrak{n}_i^{ \frac{ p_- }{ p^{\ast}_+ } } \le C \sum_{ i \in \mathfrak{I} } \mathfrak{m}_i < \infty,
   \end{equation}
the inequalities \eqref{eta3}- \eqref{n_i sum} give $ \tilde{\mathfrak{I}}= \left\{ i \in \mathfrak{I}; \, \mathfrak{n}_i > 0 \right\}$ is a finite set.  From this, one of the two possibilities below occurs: \\

  \noindent $a)$ There exists $ \mathfrak{n}_{i_{1}}, \ldots, {\mathfrak n}_{i_{s}} > 0 $ for a maximal $ s \in \mathbb{N} $; \\
  \mbox{}\\
  \noindent $b)$ $ \mathfrak{n}_i = 0 $, for all $ i \in \mathfrak{I}  $. \\

   We begin analyzing $a)$. For this, fix $ 0 < \epsilon_0 < 1 $ sufficiently small such that
   $$
      B_{ \epsilon_0 }(x_1), \cdots, B_{ \epsilon_0 }(x_s) \subset B_{ \frac{1}{\epsilon_0} }(0) \ \text{and} \ B_{ \epsilon_0 }(x_i) \cap B_{ \epsilon_0 }(x_j) = \emptyset, \ i \ne j,
   $$
  where $ x_1, \ldots, x_s $ are the singular points related to $ \mathfrak{n}_{i_1}, \ldots, \mathfrak{n}_{i_s}$, respectively. Setting 
   $$
      \psi_{ \epsilon }(x) = \varphi( \epsilon x ) - \sum_{ i=1 }^s \varphi \left( \frac{x-x_i}{\epsilon} \right) \, \forall x \in \mathbb{R}^{N},
   $$
 we have for $ 0 < \epsilon < \frac{1}{2} \epsilon_0 $, 
   $$
      \psi_{ \epsilon}(x) =
      \begin{cases}
         0, & \ \text{if} \ x \in \displaystyle \bigcup_{ i=1 }^s B_{ \frac{\epsilon}{2} }( x_i ) \\
         1, & \ \text{if} \ x \in A_{ \epsilon } = B_{ \frac{1}{\epsilon} }( 0 ) \setminus \displaystyle \bigcup_{ i=1 }^s B_{ 2 \epsilon }( x_i ),
      \end{cases}
   $$
from where it follows that  
   $$
      \text{supp} \, \psi_{ \epsilon } \subset \overline{B_{ \frac{2}{\epsilon} }( 0 )} \setminus \bigcup_{ i=1 }^s B_{ \frac{\epsilon}{2} }( x_i )
   $$   
   and  
   $$
      \int_{\mathbb{R}^{N}} | v_n |^{ p^{ \ast }(x) } \psi_{ \epsilon }  \to \int_{\mathbb{R}^{N}} | v |^{ p^{ \ast }(x) } \psi_{ \epsilon }.
   $$
   Since
   $$
      I'_{ \infty }( v_n )( v_n \psi_{ \epsilon } ) = o_n(1) \ \text{and} \ I'_{ \infty }( v_n )( v \psi_{ \epsilon } ) = o_n(1),
   $$
   repeating the same type of arguments for the case where the exponents are constant, we obtain 
   $$
      \lim_n \int_{ A_{ \epsilon } } (P_n + V(x)Q_n)=0,
   $$
   where 
 $$
 P_n( x ) = \left( \left| \nabla v_n \right|^{ p(x)-2 } \nabla v_n - \left| \nabla v \right|^{ p(x)-2 } \nabla v \right)\Big( \nabla v_n - \nabla v \Big) \,\,\, \forall x \in \mathbb{R}^{N} \,\,\, \mbox{and} \,\,\, \forall n \in \mathbb{N}. 
  $$  
and
$$
Q_n(x)=\left( \left| v_n \right|^{ p(x)-2 } v_n - \left| v \right|^{ p(x)-2 } v \right)\Big( v_n - v \Big) \,\,\, \forall x \in \mathbb{R}^{N} \,\,\, \mbox{and} \,\,\, \forall n \in \mathbb{N}.
$$
Once that  
   \begin{equation} 
      P_n(x) \ge \label{P_n inequality} 
      \begin{cases} 
         \frac{2^{ 3-p_+ }}{p_+} \left| \nabla v_n - \nabla v \right|^{ p(x) }, \ \text{if} \ p(x) \ge 2 \\
         \left( p_- - 1 \right) \frac{\left| \nabla v_n - \nabla v \right|^2}{{\left( \left| \nabla v_n \right| + \left| \nabla v \right| \right)}^{ 2-p(x) }}, \ \text{if} \ 1< p(x) < 2,
      \end{cases}
   \end{equation}
we see that  
   $$
      \int_{ A_{ \epsilon } } P_n \, dx \ge C \int_{ A_{ \epsilon } \cap \left\{ x \in \mathbb{R}^{N}; \, p(x) \ge 2 \right\} } \left| \nabla v_n - \nabla v \right|^{ p(x) } \ge 0.
   $$
 Thus, 
  \begin{equation} \label{P1}
      \lim_n \int_{  A_{ \epsilon } \cap \left\{ x \in \mathbb{R}^{N}; \, p(x) \ge 2 \right\} } \left| \nabla v_n - \nabla v \right|^{ p(x) } = 0.
  \end{equation} 
On the other hand, by H\"older's inequality
   \begin{align*}
      & \int_{  A_{ \epsilon } \cap \left\{ x \in \mathbb{R}^{N}; \, 1 <  p(x) < 2 \right\} } \left| \nabla v_n - \nabla v \right|^{ p(x) } \\
      & \phantom{ \int_{\mathbb{R}^{N}} } \le C \left| \frac{\left| \nabla v_n - \nabla v \right|^{ p(x) }}{\big( \left| \nabla v_n \right| + \left| \nabla v \right| \big)^{ \frac{p(x)(2-p(x))}{2} } } \right|_{ L^{ \frac{2}{p(x)} }  \left( \tilde{A_{ \epsilon }} \right)  }
         \! \! \! \left| \big( \left| \nabla v_n \right| + \left| \nabla v \right| \big)^{ \frac{p(x)(2-p(x))}{2} } \right|_{ L^{ \frac{2}{2-p(x)} } \left( \tilde{A_{ \epsilon }} \right) },   
  \end{align*}
   where $ \tilde{A_{ \epsilon }} = A_{ \epsilon } \cap \left\{ x \in \mathbb{R}^{N}; \, 1 <  p(x) < 2 \right\} $. From relation \eqref{P_n inequality}, the right side of above inequality goes to zero. Hence, 
\begin{equation} \label{P2}
\lim_{n}\int_{  A_{ \epsilon } \cap \left\{ x \in \mathbb{R}^{N}; \, 1 <  p(x) < 2 \right\} } \left| \nabla v_n - \nabla v \right|^{ p(x) }=0.   \\
\end{equation}  
Now,  (\ref{P1}) combined with (\ref{P2}) gives 
$$
\lim_{n}\int_{  A_{ \epsilon } } \left| \nabla v_n - \nabla v \right|^{ p(x) }=0.
$$  
The same arguments can be used to prove that
$$
\lim_{n}\int_{  A_{ \epsilon } } V(x)\left| v_n - v \right|^{ p(x) }=0.
$$ 
Therefore,
$$
v_n \to v \,\,\, \mbox{in} \,\,\, W^{1,p(x)}(A_{ \epsilon }).
$$
The last limit yields, up to a subsequence,
   $$
       \nabla v_n( x ) \to \nabla v( x ) \ \,\,\, \mbox{a.e in} \,\,\,  \ A_{\epsilon} \,\,\, ( 0 < \epsilon < \frac{1}{2} \epsilon_0).
   $$ 
Observing that  
   $$ 
      \mathbb{R}^{N} \setminus \left\{ x_1, x_2, \ldots, x_s \right\} = \bigcup_{ n \in \mathbb{N} \atop \frac{1}{n} < \frac{1}{2} \epsilon_0 } A_{ \frac{1}{n} },
   $$ 
we conclude  by a diagonal argument, that there is a subsequence of $(v_n)$, still denoted by itself, such that 
   $$
      \nabla v_n( x ) \to \nabla v( x ) \,\,\,  \mbox{a.e in} \,\, \, \mathbb{R}^{N}.
   $$
For the case $b)$, we consider
   $$
      \psi_{ \epsilon }( x ) = \varphi( \epsilon x ) \,\,\, \forall x \in \mathbb{R}^{N} \,\,\,  \text{and} \,\,  A_{ \epsilon } = B_{ \frac{1}{\epsilon} }(0), \ \epsilon > 0.
   $$
Repeating the same arguments used in the case $a)$, we have that 
\begin{equation} \label{limite}
v_n \to v \,\,\, \mbox{in} \,\,\, W^{1,p(x)}(B_{ \frac{1}{\epsilon} }(0)), \,\,\,\, \forall \epsilon >0.
\end{equation}
This way, there is again a subsequence of $(v_n)$, still denoted by itself, such that  
$$
      \nabla v_n( x ) \to \nabla v( x ) \ \,\,\,  \mbox{a.e in} \,\, \, \ \mathbb{R}^{N}.
$$
Furthermore, from \eqref{limite},  
$$
v_n \to v \,\,\, \mbox{in} \,\,\, W_{loc}^{1,p(x)}(\mathbb{R}^{N}).
$$
.\qed

\begin{corollary} \label{ESTBAIXO}  
   Let $(v_n)$ be a $(PS)_d$ sequence for $I_\infty$ with $$ d < \beta=\theta \alpha, $$ where $\theta$ and $\alpha$ were given in \eqref{theta} and \eqref{eta3} respectively. If $v_n \rightharpoonup v$ in $W^{1,p(x)}(\mathbb{R}^{N})$, then
$$
v_n \to v \,\,\, \mbox{in} \,\,\, W_{loc}^{1,p(x)}(\mathbb{R}^N).
$$ 
\end{corollary}
\noindent {\bf Proof.} Once that $d \in (-\infty, \beta)$, we claim that 
$$
\tilde{\mathfrak{I}}= \left\{ i \in \mathfrak{I}; \, \mathfrak{n}_i > 0 \right\}=\emptyset.
$$
In fact, arguing by contradiction that $\tilde{\mathfrak{I}} \not= \emptyset$, there exists $i \in \mathfrak{I}$ such that
$$
\mathfrak{n}_i \geq \alpha, 
$$
where $\alpha$ was given in \eqref{eta3}. Using the fact that $(v_n)$ is bounded in $W^{1,p(x)}(\mathbb{R}^{N})$, we have that $I'(v_n)v_n=o_n(1)$. Then, 
$$
d + o_n(1)=I_\infty(v_n)-\frac{1}{p_+}I'_\infty(v_n)v_n \geq \theta\int_{\mathbb{R}^{N}}|v_n|^{p^{*}(x)}.
$$
Letting the limit of $n \to +\infty$ in the last inequality  and using \eqref{MED2}, we get
$$
d \geq \theta \mathfrak{n}_i \geq \theta \alpha = \beta,
$$
which is an absurd. \qed

\begin{corollary} \label{NOVOPASSO} There is $\tilde{\mu}>0$ such that if $\mu \geq \tilde{\mu}$, then the $(PS)_{c_{\infty}}$ sequence $(u_n)$ given in \eqref{NIVELZ} verifies 
$$
u_n \to u \,\,\, \mbox{in} \,\,\, W_{loc}^{1,p(x)}(\mathbb{R}^{N}),
$$
where $u$ is the weak limit of $(u_n)$ in $ W^{1,p(x)}(\mathbb{R}^{N})$.
\end{corollary}

\noindent {\bf Proof.}  Using \eqref{LIMITEDONIVEL}, there is $\tilde{\mu}>0$ such that
$$
c_\infty < \beta, \,\,\, \forall \mu \geq \tilde{\mu}.
$$
Now, the corollary follows applying the Corollary \ref{ESTBAIXO}. \qed

\begin{lemma} \label{modular inequality}
   Let $ ( v_n ) $ be  a $ (PS)_d $ sequence for $ I_{ \infty } $ with
   $$ 
      d < \frac{1}{2K^{ p_+ }} \nu.
   $$
   Then, there exists $ n_0 \in \mathbb{N} $ such that
   $$
      \left( \int_{\mathbb{R}^{N}} | v_n |^{ p^{ \ast }(x) } \right)^{ 1/p^{ \ast }_- } \le K \left(  \int_{\mathbb{R}^{N}} (| \nabla v_n  |^{ p(x) } + V(x) | v_n |^{ p(x) } )\right)^{ 1/p_+ }  \,\,\, \forall n \geq n_0.
   $$ 
 
\end{lemma}

\noindent {\bf Proof.} \, Basically, we have to prove that
\begin{equation} \label{EST1}
     | v_n |_{ p^{ \ast }(x) }, \ \| v_n \| \le 1, \ \forall n \ge n_0.
\end{equation}
  If the above inequality holds,  
   $$
     \int_{\mathbb{R}^{N}} | v_n |^{ p^{ \ast }(x) }  \leq | v_n |_{ p^{ \ast }(x) }^{ p^{ \ast }_- } \,\,\,  \forall n \ge n_0
    $$
and 
$$
\| v_n \|^{ p_+ } \le \int_{\mathbb{R}^{N}} (| \nabla v_n |^{ p(x) } + V(x) | v_n |^{ p(x) })  \,\,\,  \forall n \ge n_0. 
$$
Thus, for all $n \geq n_0$
   $$
     \left( \int_{\mathbb{R}^{N}} | v_n |^{ p^{ \ast }(x) } \right)^{ 1/p^{ \ast }_- } \le | v_n |_{ p^{ \ast }(x) } \le K \| v_n \| \le K \left(  \int_{\mathbb{R}^{N}} | \nabla v_n |^{ p(x) } + V(x) | v_n |^{ p(x) } \right)^{ 1/{p_+} }.
   $$
Now, we will show that (\ref{EST1}) holds. To this end, we begin recalling  that
   $$ 
      I_{ \infty }( v_n ) - \frac{1}{q_-}I'_{ \infty }( v_n )v_n = d + o_n(1).
   $$
Therefore,
$$ 
d  + o_n(1) \geq \nu \int_{\mathbb{R}^{N}} (| \nabla v_n |^{ p(x) } + V(x) | v_n |^{ p(x) })  
$$
leading to  
   $$
      \frac{1}{\nu} d + o_n(1) \ge \int_{\mathbb{R}^{N}}( | \nabla v_n |^{ p(x) } + V(x) | v_n |^{ p(x) }).
   $$
   By the hypothesis on $ d $,
   $$ 
      \limsup_n \int_{\mathbb{R}^{N}}( | \nabla v_n |^{ p(x) } + V(x) | v_n |^{ p(x) }) \leq \frac{1}{2K^{ p_+ }}.
   $$
   Since 
   $$ 
      \frac{1}{2K^{p_+ }} < \frac{1}{K^{p_+ }},
   $$
there   exists $ n_0 \in \mathbb{N} $ such that 
   $$ 
      \int_{\mathbb{R}^{N}} (| \nabla v_n |^{ p(x) } + V(x) | v_n |^{ p(x) } ) \leq \frac{1}{K^{ p_+ }} \le 1 \,\,\, \forall n \ge n_0,
   $$
 from where it follows that  
   $$
      \| v_n \| \le 1 \,\,\, \forall  n \ge n_0. 
   $$
Thereby, 
   $$ 
      \| v_n \|^{ p_+ } \le \int_{\mathbb{R}^{N}} (| \nabla v_n |^{ p(x) } + V(x) | v_n |^{ p(x) } ) \,\,\,\, \forall  n \ge n_0,
   $$
implying that  
$$
      \| v_n \| \leq \frac{1}{K}<1 \,\,\,\, \forall  n \ge n_0.
$$      
Once that   
$$     
| v_n |_{ p^{ \ast }(x) } \leq K \| v_n \|  \,\,\,\, \forall n,
$$
it follows that
$$
| v_n |_{ p^{ \ast }(x) } \leq 1 \,\,\,\, \forall  n \ge n_0.
$$ \qed

\begin{lemma} \label{(PS)_d behaviour}
   Let $ ( v_n ) $ be  a $ (PS)_d $ sequence for $ I_{ \infty } $ with
   $$ 
      d < \min \left \{ \theta \left( \frac{1}{K} \right)^{ \frac{1}{\theta} }, \frac{1}{2K^{ p_+ }} \nu \right \}.
   $$
   Then, up to a subsequence, \\
   
   \noindent $(a)$ \,\,\, $ v_n \to 0 $ in $ W^{ 1,p(x) } \big( \mathbb{R}^{N} \big) $, or\\
   \mbox{}\\
   \noindent $(b)$ There are $ R $, $ \eta > 0 $ and $ (y_n) \subset \mathbb{R}^{N} $ such that
            $$
              \limsup_n \int_{ B_R( y_n ) } | v_n |^{ p(x) } \geq \eta.
            $$ 

\end{lemma}

\noindent {\bf Proof.} Up to a subsequence, we can assume that
   $$ 
      \int_{\mathbb{R}^{N}} (| \nabla v_n |^{ p(x) } + V(x) | v_n |^{ p(x) }) \to L \ge 0.
   $$
   If $ ( b ) $ does not hold, there is $ R > 0 $ such that 
   $$ 
      \lim_{n} \sup_{ y \in \mathbb{R}^{N} } \int_{ B_R(y) } | v_n |^{ p(x) } = 0.
   $$ 
   Since $ ( v_n ) $ is bounded, by Lemma $ 3.1 $ in \cite{FZZ}, 
   $$ 
      v_n \to 0 \ \text{in} \ L^{ q(x) } \big( \mathbb{R}^{N} \big)
   $$
 or equivalently
\begin{equation}\label{L1} 
      \int_{\mathbb{R}^{N}} | v_n |^{ q(x) }  \to 0.
\end{equation}
   Our goal is to prove that $ L = 0 $, because if this occurs, we get $ ( a ) $.
   Suppose by contradiction that $ L > 0 $. Since $ I'_{ \infty }( v_n )v_n = o_n(1) $, we derive that
   $$ 
      \int_{\mathbb{R}^{N}} | v_n |^{ p^{ \ast }(x) }  \to L.
   $$
By (\ref{L1}), 
   \begin{align*}
       d + o_n(1) & = I_{ \infty }( v_n ) + \mu \int_{\mathbb{R}^{N}} \frac{1}{q(x)} | v_n |^{ q(x) }  \\
      & = \int_{\mathbb{R}^{N}} \frac{1}{p(x)} \left( | \nabla v_n |^{ p(x) } + V(x) | v_n |^{ p(x) } \right)  - \int_{\mathbb{R}^{N}} \frac{1}{p^{ \ast }(x)} | v_n |^{ p^{ \ast }(x) }
   \end{align*}
 and so,
   $$
      d + o_n(1) \ge \frac{1}{p_+} \int_{\mathbb{R}^{N}} (| \nabla v_n |^{ p(x) } + V(x) | v_n |^{ p(x) } ) - \frac{1}{p^{ \ast }_-} \int_{\mathbb{R}^{N}} | v_n |^{ p^{ \ast }(x) } .
   $$
Taking the limit of $n \to +\infty$ in the last inequality, 
   \begin{equation} \label{relation3}
      d \ge \frac{1}{p_+} L - \frac{1}{p^{ \ast }_-} L = \theta L.
   \end{equation}
   On the other hand, using the fact that  
   $$ 
      d < \frac{1}{2K^{ p_+ }} \nu,
   $$
it follows from  Lemma \ref{modular inequality} that there exists $ n_0 \in \mathbb{N} $ such that
   $$
      \left( \int_{\mathbb{R}^{N}} | v_n |^{ p^{ \ast }(x) }  \right)^{ 1/p^{ \ast }_- } \le K \left(  \int_{\mathbb{R}^{N}} | \nabla v_n |^{ p(x) } + V(x) | v_n |^{ p(x) } \right)^{ 1/p_+ }, \,\,\, \forall n \ge n_0.
   $$ 
   Taking the limit of $n \to +\infty$, we derive 
   $$
      L^{ 1/p^{ \ast }_- } \le K L^{ 1/p_+}
   $$
or equivalently  
   \begin{equation} \label{relation4}
      L \ge \left( \frac{1}{K} \right)^{ \frac{1}{\theta} }.
   \end{equation}
Combining \eqref{relation3} with \eqref{relation4}, 
   $$
      d \ge \theta L \ge \theta \left( \frac{1}{K} \right)^{ \frac{1}{\theta} },
   $$
   which is a contradiction, showing that $ L = 0 $. \qed

\vspace{1 cm}

\subsection{Proof of Theorem \ref{T2}}

Consider the $ (PS)_{ c_{ \infty } } $ sequence $ ( u_n ) $ for $ I_{ \infty } $ given in  \eqref{NIVELZ} and $ u \in W^{ 1,p(x) } \big( \mathbb{R}^{N} \big) $ such that $ u_n \rightharpoonup u $ in $W^{ 1,p(x) } \big( \mathbb{R}^{N} \big) $. Provided that $ u \ne 0 $, by Lemma \ref{aew convergence}, we have a solution for $ (P_\infty) $. Now, if $ u = 0 $, we need to work a little more. We recall that there is $ \mu \ge \mu_{ \infty } $ such that 
   $$ 
       0<c_{ \infty } < \min \left \{ \theta \left( \frac{1}{K} \right)^{ \frac{1}{\theta} }, \frac{1}{2K^{ p_+ }} \nu \right \}.
   $$
By Lemma \ref{(PS)_d behaviour}, we know that $ (b) $ holds, that is, there exist $ R $, $ \eta > 0 $  and $ ( y_n ) \subset \mathbb{R}^{N} $ such that
   $$
      \limsup_n \int_{ B_R( y_n ) } | u_n |^{ p(x) } \geq \eta.
   $$ 
Hereafter, without lost of generality, we assume that $(y_n) \subset \mathbb{Z}^{N}$ and define   
   $$
      \widetilde{u_n}(x) = u_n( x+y_n ) \,\,\, \forall   x \in \mathbb{R}^{N}.
   $$
A straightforward  calculus gives 
   $$
      \ I_{ \infty }( \widetilde{u_n} ) = I_{ \infty }( u_n ) \ \text{and} \ I'_{ \infty }( \widetilde{u_n} ) \to 0,
   $$
 showing that $ ( \widetilde{u_n} ) $ is also a $ (PS)_{ c_{ \infty } } $ sequence for $ I_{ \infty } $. In what follows, we denote by $ \widetilde{u} \in W^{ 1,p(x) }( \mathbb{R}^{N} ) $ the weak limit of $ (\widetilde{u_n} )$.  Once that 
   $$
      \int_{ B_R(0) } | \widetilde{u_n} |^{ p(x) } = \int_{ B_R( y_n ) } | u_n |^{ p(x) } 
   $$
   and 
   $$
      W^{ 1,p(x) } \big( \mathbb{R}^{N} \big) \overset{compactly}{\hookrightarrow} L^{ p(x) } \big( B_R(0) \big), 
   $$
   it follows that
   $$ 
      \int_{ B_R(0) } | \widetilde{u} |^{ p(x) } \ge \eta > 0,
   $$
   implying that $ \widetilde{u} \ne 0 $.  Thereby, by Lemma \ref{aew convergence}, $ \widetilde{u} $ is a solution for $ ( P_{ \infty } ). $ 

\vspace{.7cm}

In the sequel, $u_{ \infty } \in W^{1,p(x)}(\mathbb{R}^{N})$ denotes the solution found by the above arguments. Moreover, let us denote by $ {\mathcal N}_{ \infty } $ the Nehari manifold  associated with $ I_{ \infty } $ given by 
   $$
      {\mathcal N}_{ \infty } = \left\{ u \in W^{ 1,p(x) } \big( \mathbb{R}^{N} \big) \setminus \{ 0 \}; \, I'_{ \infty }( u )u = 0 \right\}. \vspace{.2cm}
   $$
Our goal is to show that $u_{ \infty }$ is a {\it ground state solution}, that is,  
   $$ 
      I_{ \infty }(u_{ \infty }) = \inf_{ u \in {\mathcal N}_{ \infty } } I_{ \infty }(u). 
   $$
Moreover, we will prove that $u_{ \infty }$ can be chosen as a nonnegative solution.  

\begin{proposition} For all $u \in \mathcal{N}_{\infty}$, we have that $I_{\infty}(u)>0$ and 
     $$
      0<J_{ \infty } = \inf_{ u \in \mathcal{N}_{ \infty } } I_{ \infty }(u).
   $$
  
\end{proposition}

\noindent {\bf Proof.}   If $ u \in \mathcal{N}_{ \infty } $,
   
   $$
       \int_{\mathbb{R}^{N}} (| \nabla u |^{ p(x) } + V(x) | u |^{ p(x) }) = \mu \int_{\mathbb{R}^{N}} | u |^{ q(x) } + \int_{\mathbb{R}^{N}} | u |^{ p^{ \ast }(x) }.
   $$ 
 This way, 
$$
      I_{ \infty }(u) \geq  \mu \nu \int_{\mathbb{R}^{N}} | u |^{ q(x) }  + \theta \int_{\mathbb{R}^{N}} | u |^{ p^{ \ast }(x) }  > 0.
 $$
Now, arguing by contradiction, if $ J_{ \infty } = 0 $, there exists $ ( v_n ) \subset \mathcal{N}_{ \infty } $ such that
\begin{equation} \label{N0}
      I_{ \infty }( v_n ) \to 0.
 \end{equation}
   Since $ v_n \in \mathcal{N}_{ \infty },$ by the last inequality  
    \begin{equation} \label{N2}
      I_{ \infty }( v_n)  \geq \mu \nu \int_{\mathbb{R}^{N}} | v_n |^{ q(x) }  + \theta \int_{\mathbb{R}^{N}} | v_n |^{ p^{ \ast }(x) } .
   \end{equation}
 Using (\ref{N0}) together with \eqref{N2}, we derive the limits 
   $$
     \int_{\mathbb{R}^{N}} | v_n |^{ q(x) } , \; \int_{\mathbb{R}^{N}} | v_n |^{ p^{ \ast }(x) }  \to 0
   $$
which lead to  
   $$
      \int_{\mathbb{R}^{N}}( | \nabla v_n |^{ p(x) } + V(x) | v_n |^{ p(x) }) \to 0,
   $$
   or equivalently,
   $$
      \| v_n \| \to 0.
   $$
  On the other hand, by Sobolev embedding, there are constants $ C_1, C_2 > 0 $ verifying 
$$
      | u |_{ q(x) } \le C_1 \| u \| \,\,\, \mbox{and} \,\,\,  | u |_{ p^{ \ast }(x) } \le C_2 \| u \| \,\,\, \forall u \in W^{1,p(x)}(\mathbb{R}^N).
$$
Consequently, there is $n_0 \in \mathbb{N}$ such that 
   $$
     \| v_n \|, | v_n |_{ q(x) }, | v_n |_{ p^{ \ast }(x) } \le 1 \,\,\,\, \forall n \geq n_0 .
   $$
The last inequalities together with the fact that $(v_n) \subset  \mathcal{N}_\infty$ give
$$
     \| v_n \|^{ p_+ } \leq \mu | v_n |_{ q(x) }^{ q_- } + | v_n |_{ p^{ \ast}(x) }^{ p^{ \ast }_- } \le \mu C_3 \| v_n \|^{ q_- } + C_4 \| v_n \|^{ p^{ \ast }_- }, \, \forall n \ge n_0, 
$$ 
for some positive constants $C_3$ and $C_4$. Thus,  
   $$
      1 \le \mu C_3 \| v_n \|^{ q_- - p_+ } + C_4 \| v_n \|^{ p^{ \ast }_- - p_+ }, \, \forall n \ge n_0,
   $$
obtaining a contradiction, because $ \| v_n \| \to 0 $. Then, $J_{ \infty } > 0.$ \qed

\begin{corollary}
   Any ground state solution $ u $ to $ ( P_{ \infty } ) $ has a well defined sign, that is, $ u \ge 0 $ or $ u \le 0 $.    
\end{corollary}

\noindent {\bf Proof.} Using the fact that $u$ is a solution, if follows that  
$$
I_\infty '(u^\pm)u^\pm=0,
$$
where $u^{+}=\max\{u,0\}$ and $u^{-}=\min\{u,0\}$. Thus, if $u^{\pm}\not=0$, we see that $u^{\pm} \in \mathcal{N}_{ \infty }$, and so, 
$$
J_\infty=I_\infty(u)= I_\infty(u^+)+I_\infty(u^-) \geq 2 J_\infty, 
$$
which is an absurd. Therefore, $u^-=0$ or $u^+=0$. \qed

\begin{proposition} \label{intersection property}
   Let $ u \in W^{ 1,p(x) } \big( \mathbb{R}^{N} \big) \setminus \{ 0 \} $. Then, exists a unique $ t_u > 0 $ such that 
   $$ 
      t_u u \in \mathcal{N}_{ \infty }.
   $$ 
\end{proposition}

\noindent {\bf Proof.}  The proof follows the same arguments explored in \cite{Fan}. \qed

\begin{proposition} \label{Nehari Infimum}
   The mountain pass level $ c_{ \infty } $ satisfies
   $$
      c_{ \infty } = J_{ \infty } = \inf_{ u \in \mathcal{N}_{ \infty } } I_{ \infty }(u).
   $$
\end{proposition}

\noindent {\bf Proof.}   Let $ u \in \mathcal{N}_{ \infty } $ and choose $ t_0 > 0 $ such that $ u_0 = t_0 u $ satisfies $ I_{ \infty }( u_0 ) < 0 $. Then
   $$
      \gamma_0(t) = t u_0, \,\,\, \forall  t \in [0,1],
   $$
   belongs to $ \Gamma $. Hence, 
   $$ 
      c_{ \infty } \le \max_{ t \in [0,1] } I_{ \infty }( \gamma_0(t) ) = \max_{ s \in [0,t_0] } I_{ \infty }(su) \le \max_{ s \ge 0 } I_{ \infty }(su) = I_{ \infty }(u). 
   $$
   Thereby,
 \begin{equation} \label{cinfinito1}
      c_{ \infty } \le J_{ \infty }.
\end{equation}
 
   For the reversed inequality, consider a $ (PS)_{ c_{ \infty } } $ sequence $ ( u_n ) $ for $ I_{ \infty } $. Since $ c_{ \infty } > 0 $, we can assume that 
   $ u_n \ne 0 $ for all $n \in \mathbb{N}$. In this case, by Proposition \ref{intersection property}, for each $ n \in \mathbb{N} $, there exists a unique
    $ t_n > 0 $ such that $ t_n u_n \in \mathcal{N}_{ \infty }.$ Thus, 
   \begin{equation} \label{N4}
      \int_{\mathbb{R}^{N}} t_n^{ p(x) } \left( | \nabla u_n |^{ p(x) } + V(x) | u_n |^{ p(x) } \right) 
      = \mu \int_{\mathbb{R}^{N}} t_n^{ q(x) } | u_n |^{ q(x) } + \int_{\mathbb{R}^{N}} t_n^{ p^{ \ast }(x) } | u_n |^{ p^{ \ast }(x) } .
   \end{equation}
   From \eqref{N4}, we see that $ t_n \not \to 0. $ In fact, if $ t_n \to 0 $, we can assume $ t_n < 1, \, \forall n \in \mathbb{N} $, and so, from \eqref{N4},
   $$
      t_n^{ p_+ } \int_{\mathbb{R}^{N}}( | \nabla u_n |^{ p(x) } + V(x) | u_n |^{ p(x) }) \le \mu t_n^{ q_- } \int_{\mathbb{R}^{N}}  | u_n |^{ q(x) }  + t_n^{ p^{ \ast }_- } \int_{\mathbb{R}^{N}} | u_n |^{ p^{ \ast }(x) } .
   $$
   Since $ t_n > 0$ for all $n \in \mathbb{N} $, we get
   \begin{equation} \label{N5}
      \int_{\mathbb{R}^{N}} (| \nabla u_n |^{ p(x) } + V(x) | u_n |^{ p(x) }) \\
      \le \mu t_n^{ q_- - p_+ } \int_{\mathbb{R}^{N}}  | u_n |^{ q(x) }  + t_n^{ p^{ \ast }_- - p_+ } \int_{\mathbb{R}^{N}} | u_n |^{ p^{ \ast }(x) } .
   \end{equation}
The boundedness of $ ( u_n ) $ in $W^{1,p(x)}(\mathbb{R}^{N})$ together with Sobolev embedding yields $(u_n)$ is bounded in  $L^{ q(x) } ( \mathbb{R}^{N})$ and  $L^{ p^{ \ast }(x) }(\mathbb{R}^{N})$. Therefore, from  \eqref{N5}
   $$
      u_n \to  0 \,\,\, \mbox{in} \,\,\, W^{1,p(x)}(\mathbb{R}^{N}),
   $$
leading to $I(u_n) \to I(0)=0$, which is a contradiction, because \linebreak $ I(u_n) \to c_{ \infty } > 0 $. Analogously, from \eqref{N4}, we conclude $(t_n)$ is bounded. Indeed, if there exists a subsequence of $ ( t_n ) $, still denote by itself, satisfying $ t_n \to \infty $, we can assume $ t_n > 1$ for all $n \in \mathbb{N} $, and so, from \eqref{N4},
$$
      \int_{\mathbb{R}^{N}} (| \nabla u_n |^{ p(x) } + V(x) | u_n |^{ p(x) } ) \geq t_n^{ p^{ \ast }_- - p_+ } \int_{\mathbb{R}^{N}} | u_n |^{ p^{ \ast }(x) }.
$$
   Hence
   $$
     \frac{1}{t_n^{ p^{ \ast }_- - p_+ }} \int_{\mathbb{R}^{N}} | \nabla u_n |^{ p(x) } + V(x) | u_n |^{ p(x) }  \geq \int_{\mathbb{R}^{N}} | u_n |^{ p^{ \ast }(x) } ,
   $$
   and so,
\begin{equation} \label{E1}
      \int_{\mathbb{R}^{N}} | u_n |^{ p^{ \ast }(x) } \, dx \to 0.
\end{equation}
   By interpolation \cite[Lemma 4.1]{Alves1} , 
\begin{equation}\label{E2} 
      \int_{\mathbb{R}^{N}} | u_n |^{ q(x) } \to 0.
\end{equation}
   Since 
   \begin{equation} \label{N6}
      \int_{\mathbb{R}^{N}}( | \nabla u_n |^{ p(x) } + V(x) | u_n |^{ p(x) }) = \mu  \int_{\mathbb{R}^{N}} | u_n |^{ q(x) }  + \int_{\mathbb{R}^{N}} | u_n |^{ p^{ \ast }(x) } + o_n(1),
   \end{equation}
the limits \eqref{E1} and \eqref{E2} combined with \eqref{N6} give 
   $$
      u_n \to 0 \,\,\, \mbox{in} \,\,\, W^{1,p(x)}(\mathbb{R}^{N}),
   $$
implying that $I_\infty(u_n) \to0$,  which is a contradiction, because \linebreak $I_\infty(u_n) \to c_\infty >0$. Then,  $ (t_n ) $ is bounded and, up to a subsequence, there exists
$ t_0 \in (0, \infty) $ such that $ t_n \to t_0 $. Our goal is to prove that $ t_0 = 1 $,  because if it is  true, we get
   \begin{align*}
      J_{ \infty } & \le I_{ \infty }( t_n u_n ) \\
      & = \int_{\mathbb{R}^{N}} \frac{1}{p(x)} \left( | \nabla t_n u_n |^{ p(x) } + V(x) | t_n u_n |^{ p(x) } \right) \ - \mu  \int_{\mathbb{R}^{N}} \frac{1}{q(x)} | t_n u_n |^{ q(x) } \\ 
      & \phantom{ = } - \int_{\mathbb{R}^{N}} \frac{1}{p^{ \ast }(x)} | t_n u_n |^{ p^{ \ast }(x) } \\
      & \le I_{ \infty } (u_n) + a(n) \int_{\mathbb{R}^{N}} \frac{1}{p(x)} \left( | \nabla u_n |^{ p(x) } + V(x) | u_n |^{ p(x) } \right) \\
      & \phantom{ \le } - \mu b(n) \int_{\mathbb{R}^{N}} \frac{1}{q(x)} | u_n |^{ q(x) }  - c(n) \int_{\mathbb{R}^{N}} \frac{1}{p^{ \ast }(x)} | u_n |^{ p^{ \ast }(x) } \\\
      & =  c_{ \infty } + o_n(1),  
   \end{align*}
   where
   $$
      a(n) = \max \left\{ t_n^{ p_+ },t_n^{ p_- } \right\} -1, \quad b(n) = \min \left\{ t_n^{ q_+ },t_n^{ q_- } \right\} -1,  
   $$ 
and 
$$
    c(n) = \min \left\{ t_n^{ p^{ \ast }_+ },t_n^{ p^{ \ast }_- } \right\} - 1 .  
$$   
   Hence, taking the limit of $n \to +\infty,$ we obtain
\begin{equation} \label{cinfinito2}
      J_{ \infty } \le c_{ \infty }.
\end{equation}
From \eqref{cinfinito1} and \eqref{cinfinito2} , it follows that $c_\infty = J_\infty$. 

In what follows, we will show that $t_0=1$. In fact, if $ t_0 > 1 $, we can assume without lost of generality that $ t_n > 1, \, \forall \, n \in \mathbb{N} $. Thus, by \eqref{N4}, 
   \begin{equation} \label{N7}
      \int_{\mathbb{R}^{N}}( | \nabla u_n |^{ p(x) } + V(x) | u_n |^{ p(x) }) \geq \mu t_n^{ q_- - p_+ } \int_{\mathbb{R}^{N}} | u_n |^{ q(x) }+ t_n^{ p^{ \ast }_- - p_+ } \int_{\mathbb{R}^{N}} | u_n |^{ p^{ \ast }(x) }.
   \end{equation}
From \eqref{N6} and \eqref{N7}, 
   $$
      0 \ge \mu \left( t_n^{ q_- - p_+ } - 1 \right) \int_{\mathbb{R}^{N}}  | u_n |^{ q(x) }  + \left( t_n^{ p^{ \ast }_- - p_+ } - 1 \right) \int_{\mathbb{R}^{N}} | u_n |^{ p^{ \ast }(x) }  + o_n(1).
   $$
   So,
   $$
      0 \ge \mu \left( t_0^{ q_- - p_+ } - 1 \right) L_1 + \left( t_0^{ p^{ \ast }_- - p_+ } - 1 \right) L_2,
   $$
   where 
   $$ 
      0 \le L_1 = \lim_n \int_{\mathbb{R}^{N}} | u_n |^{ q(x) }  \ \text{and} \ 0 < L_2 = \lim_n \int_{\mathbb{R}^{N}} | u_n |^{ p^{ \ast }(x) } . 
   $$
   Once that $ t_0 > 1 $, 
   $$ 
      \mu \left( t_0^{ q_- - p_+ } - 1 \right) L_1 + \left( t_0^{ p^{ \ast }_- - p_+ } - 1 \right) L_2 > 0,
   $$
   which is an absurd.  The case $ t_0 < 1 $ can be studied of the same way. Thereby, $t_0 = 1.$ \qed

As an immediate consequence of the last result, we have the corollary below
\begin{corollary}
   The mountain pass level $ c_{ \infty } $ also satisfies
   $$
      c_{ \infty } =  \inf_{ u \in W^{ 1,p(x) } ( \mathbb{R}^{N} ) \atop u \ne 0 } \max_{ t \ge 0 } I_{ \infty }(tu).
   $$
\end{corollary}

The proof of next corollary follows the ideas explored in \cite{W} and its proof we be omitted.  

\begin{corollary} \label{min is gstate}
   If $ u \in \mathcal{N}_\infty $ and $ I_\infty(u)=J_\infty $, then $ u $ is a ground state solution.
\end{corollary}

The next corollary ensures that the solution $ u_{ \infty } $ found is a ground state solution. As consequence, $  u_{ \infty } $ has a well defined sign.

\begin{corollary}
   Let $ ( u_n ) $ be  a $ (PS)_{ c_{ \infty } } $ sequence for $ I_{ \infty } $ and $ u \in W^{ 1,p(x) } \big( \mathbb{R}^{N} \big) $ such that $ u_n \rightharpoonup u $ in $ W^{ 1,p(x) } \big( \mathbb{R}^{N} \big) $. Then, if $ u \ne 0 $, we have the equality 
   $$ 
      I_{ \infty }( u ) = J_{ \infty }.
   $$ 
\end{corollary}

\noindent {\bf Proof.}  Once that we are assuming that  $ u \ne 0 $, the Lemma \ref{aew convergence} implies that $ u \in \mathcal{N}_{ \infty } $. Therefore,
   $$ 
      I_{ \infty }( u ) \ge J_{ \infty }.
   $$
   On the other hand, by Proposition \ref{Nehari Infimum}, 
   \begin{align*}
      J_{ \infty } & = c_{ \infty } = I_{ \infty }( u_n ) - \frac{1}{p_+} I'_{ \infty }( u_n )u_n + o_n(1) \\
      & = \int_{\mathbb{R}^{N}} \left( \frac{1}{p(x)} - \frac{1}{p_+} \right) \left( | \nabla u_n |^{ p(x) } + V(x) | u_n |^{ p(x) } \right) \\
      & \phantom{ = } + \mu \int_{\mathbb{R}^{N}} \left( \frac{1}{p_+} - \frac{1}{q(x)} \right) | u_n |^{ q(x) } + \int_{\mathbb{R}^{N}} \left( \frac{1}{p_+} - \frac{1}{p^{ \ast }(x)} \right) | u_n |^{ p^{ \ast }(x) }+ o_n(1) .
   \end{align*}
By Fatou's lemma, 
   \begin{multline*}
      J_{ \infty } \ge  \int_{\mathbb{R}^{N}} \left( \frac{1}{p(x)} - \frac{1}{p_+} \right) \left( | \nabla u |^{ p(x) } + V(x) | u |^{ p(x) } \right)  + \mu \int_{\mathbb{R}^{N}} \left( \frac{1}{p_+} - \frac{1}{q(x)} \right) | u |^{ q(x) }  \\
      + \int_{\mathbb{R}^{N}} \left( \frac{1}{p_+} - \frac{1}{p^{ \ast }(x)} \right) | u |^{ p^{ \ast }(x) }  = I_{ \infty }(u) - \frac{1}{p_+} I'_{ \infty }(u)u = I_{ \infty }(u), 
   \end{multline*}
from where it follows that $ I_{ \infty }( u ) = J_{ \infty }.$ \qed \\

Now, to conclude the proof of Theorem \ref{T2}, we observe that 
   $$
   I_\infty(-v)=I_\infty(v) \,\,\, \mbox{and} \,\,\, I_\infty'(-v)=-I_\infty'(v) \,\,\, \forall v \in W^{1,p(x)}(\mathbb{R}^{N}),
   $$
which ensures that $ u_{ \infty } $ can be chosen nonnegative.

\section{Proof of Theorem \ref{T1}}

In this section, we will show that the functional $I:W^{1,p(x)+\sigma(x)}(\mathbb{R}^{N}) \to \mathbb{R}$ given by 
$$
I(u)=\int_{\mathbb{R}^{N}} \frac{1}{(p(x)+\sigma(x))} \left( | \nabla u |^{ p(x)+\sigma(x) } + (V(x)-W(x)) | u |^{ p(x)+\sigma(x) } \right) - \Psi(u) - J(u) 
$$
where
$$
\Psi(u)= \mu \int_{\mathbb{R}^{N}} \frac{1}{(q(x)-\tau(x))} | u |^{ q(x)-\tau(x) } \,\,\, \mbox{and} \,\,\, J(u)=\int_{\mathbb{R}^{N}} \frac{1}{p^{ \ast }(x)} |u|^{ p^{ \ast}(x) },  
$$
has a ground state solution for $\mu$ large enough. Hereafter, we denote by $r(x), s(x)$ and $U(x)$ the functions
$$
r(x)=p(x)+\sigma(x), \,\,\, s(x)=q(x)-\tau(x) \,\,\, \mbox{and} \,\,\, U(x)=V(x)-W(x) \,\,\, \forall x \in \mathbb{R}^{N}.
$$ 
Hence, we can rewrite the functional $I$ of the following way
$$
I(u)=\int_{\mathbb{R}^{N}} \frac{1}{r(x)} \left( | \nabla u |^{r(x)} + U(x) | u |^{r(x)} \right) - \mu \int_{\mathbb{R}^{N}} \frac{1}{s(x)} | u |^{s(x)}-\int_{\mathbb{R}^{N}} \frac{1}{p^{ \ast }(x)} |u|^{ p^{ \ast}(x) } 
$$
for all $u \in W^{1,r(x)}(\mathbb{R}^{N})$. 

Using well known arguments, it follows that $ I \in C^{1}(W^{1,r(x)}(\mathbb{R}^{N}), \mathbb{R})$ with
$$
I'( u ) v  = \int_{\mathbb{R}^{N}}( | \nabla u |^{ r(x)-2 } \nabla u \nabla v + U(x) | u |^{ r(x)-2 } u v ) - \mu \int_{\mathbb{R}^{N}} | u |^{ s(x)-2 } u v - \int_{\mathbb{R}^{N}} | u |^{ p^{ \ast }(x)-2 } u v,  
$$
for all $ u,v \in W^{ 1,r(x)}( \mathbb{R}^{N} ) $. Moreover, we also have that $I$ verifies the mountain pass geometry, this way, there is $(u_n) \subset W^{1,r(x)}(\mathbb{R}^{N})$ satisfying 
$$
I(u_n) \to c \,\,\ \mbox{and} \,\,\ I'(u_n) \to 0
$$
where $c$ is the mountain pass level associated with $I$. Using the same arguments explored in the proof of Lemma   
\ref{aew convergence} and Corollary \ref{NOVOPASSO}, there is $\tilde{\mu} >0$ such that, for $\mu \geq \tilde{\mu}$,  
\begin{equation} \label{Limite1}
u_n \to u \,\, \mbox{in} \,\,\, W^{1,r(x)}_{loc}(\mathbb{R}^{N}),
\end{equation}
where $ u \in W^{ 1,r(x) }( \mathbb R^N ) $ is the weak limit of $ (u_n) $, and so, 
\begin{equation} \label{Limite2}
u_n \to u \,\, \mbox{in} \,\,\, W^{1,p(x)}_{loc}(\mathbb{R}^{N}).
\end{equation}
The limit \eqref{Limite1} leads to  
$$
I'(u)=0,
$$
showing that $u$ is a solution for $(P)$. If $u\not=0$, the same ideas used in the previous section give  $I(u)=c,$ implying that $u$ is a ground state solution for $(P)$. Now, if $u=0$, the above limits yield 
$$
I_\infty(u_n) \to c \,\,\, \mbox{and} \,\,\, I'_\infty(u_n) \to 0.
$$
The arguments used in the proof of Proposition \ref{Nehari Infimum} works to show that   
\begin{equation} \label{ESTM1}
c \geq c_\infty.
\end{equation}
Fixing $\mu \geq \max\{\mu_\infty , \tilde{\mu} \}=\mu^{*}$, we know that functional $I_\infty$ has a nonnegative ground state solution $w \in W^{1,p(x)}( \mathbb R^N )$, that is,
$$
I_\infty(w)=c_\infty \,\,\, \mbox{and} \,\,\, I'_\infty(w)=0.
$$
Moreover, for each $n \in \mathbb{N}$, we set $x_n=(n,0,...,0) \in \mathbb{R}^{N}$ and 
$$
w_n(x)=w(x+x_n) \,\,\, \forall x \in \mathbb{R}^{N}.
$$
Using Moser iteration method \cite{Moser} ( see also \cite{AlvesGio,Li} ), there is $C>0$ such that
$$
|w_n|_{L^{\infty}(B_{R_2}(z))} \leq C |w|_{L^{m^{*}}(B_{R_1}(z+x_n))} \,\,\, \mbox{for} \,\, R<R_2<R_1 \,\,\, \mbox{and} \,\, \forall n \in \mathbb{N}.
$$
Once that $w\in L^{p^{*}(x)}(\mathbb{R}^{N})$ and $p^*(x)=m^*$ in $B_{R_1}(z+x_n)$ for all $n \in \mathbb{N}$, we deduce 
$$
|w|_{L^{m^{*}}(B_{R_1}(z+x_n))} \to 0 \,\,\, \mbox{as} \,\,\, n \to +\infty.
$$ 
Then,
$$
|w_n|_{L^{\infty}(B_{R_2}(z))} \to 0 \,\,\, \mbox{as} \,\,\, n \to +\infty. 
$$
The above limit implies that 
$$
|\nabla w_n|_{L^{\infty}(B_{R_2}(z))} \to 0 \,\,\, \mbox{as} \,\,\, n \to +\infty \,\,\, ( see \, \cite{Acerbi1, Acerbi2, Fan2} ).
$$
From this, we fix $n \in \mathbb{N}$ such that
$$
|w_n(x)|,|\nabla w_n(x)| \leq 1 \,\,\, \forall x \in B_R(z).  
$$ 
Now, a simple computation gives 
$$
I'(w_n)w_n \leq I'_\infty(w_n)w_n=I'_\infty(w)w=0.
$$
Thus, there is $t_n \in (0,1]$ such that 
$$
I(t_nw_n)= \max_{t \geq 0}I(tw_n) \,\,\, \mbox{and} \,\,\, I(t_nw_n) \geq c. 
$$
Using the definition of $w_n$ together with the fact that $t_n \in (0,1]$, we have that
\begin{equation} \label{ESTM2}
c \leq I(t_nw_n) \leq I_\infty(t_nw_n) \leq I_\infty(w)=c_\infty.
\end{equation}
Combining \eqref{ESTM1} and \eqref{ESTM2}, it follows that
$$
c=I(t_nw_n).
$$
Since 
$$ 
t_nw_n \in \mathcal{N}=\{u \in W^{1,r(x)}(\mathbb{R}^{N})\setminus{0}\,;\, I'(u)u=0\}, 
$$ 
the same arguments used in the proof of Corollary \ref{min is gstate} can be used to prove that $t_nw_n$ is a ground state solution for $I$, finishing the proof of Theorem \ref{T1}. \qed

\end{document}